\documentclass{amsart}
\usepackage{mathdots}
\usepackage{bbm}
\usepackage{epsfig}
\usepackage{epstopdf}

\newtheorem{theorem}{Theorem}
\newtheorem{proposition}[theorem]{Proposition}

\newtheorem{corollary}[theorem]{Corollary}

\theoremstyle{definition}
\newtheorem{definition}[theorem]{Definition}

\usepackage{srcltx,amssymb}
\usepackage{graphicx}
\usepackage{amsthm}
\usepackage{amsmath}
\usepackage{amscd}
\usepackage{slashed}
\numberwithin{theorem}{section} 
\numberwithin{equation}{section} 

\begin{document}

    \title{A formulation of Euclidean geometry leading to recognition of a new class of algebra isomorphism invariants}
\author{Fred Greensite}

\begin{abstract}We present alternative postulates for Euclidean geometry whose merit is that they lead to a new class of invariants and associated geometries for real finite-dimensional unital associative algebras.
\end{abstract} 

\maketitle   

     \markboth{F. Greensite }{The Pythagorean Theorem and generalization of the usual algebra norm}

\section{Introduction}

The twentieth century saw multiple rigorous formulations of Euclidean geometry, three important examples of which are those of David Hilbert \cite{hilbert:1950}, Alfred Tarski \cite{tarski:1999}, and George Birkhoff \cite{birkhoff:1932}.  Each made unique contributions.  \
\begin{itemize}\item Hilbert provided a treatment of essentially the entirety of what was intended by Euclidean plane geometry. \item Tarski's formulation is less encompassing, but has the laudable features of being provably consistent, complete, and decidable. \item  Birkhoff introduced a system of (only four identified) postulates suitable to the needs of high school students (``based on scale and protractor", as he put it). The geometry is implicitly analytic given its \mbox{incorporation of $\mathbb{R}$}. \end{itemize}  

Accordingly, though without using such language, Birkhoff ultimately demonstrated that his postulates imply an affine space whose vector space component is the real Cartesian plane (though the proofs in his paper are generally synthetic in character).  In isolation, the latter entity is economically formulated from two postulates. \begin{enumerate} \item The plane is a real two-dimensional vector space with a nonnegative length function $\ell:\mathbb{R}^2\rightarrow\mathbb{R}$ which outputs the length of the line segment defined by the origin and any point of the plane. \item In the standard basis, with a point expressed by $(x,y)$, the length function is $\ell(x,y) = \sqrt{x^2+y^2}$. The length of a line segment with endpoints $P_1$, $P_2$ is $\ell(P_1-P_2)$. 
\end{enumerate}
    Of course, the seeming simplicity of this two postulate formulation of the Cartesian plane is belied by the fact that Postulate (2) is motivated by the prior existence of synthetic proofs of the Pythagorean Theorem, and the fact that Postulate (1) relies on the postulates underlying $\mathbb{R}$.  That is, all of Real Analysis is implied.
      
 In the context of these offerings by the great mathematicians Euclid, Descartes, Hilbert, Tarski, and Birkhoff, we certainly don't propose that the postulates below should be placed along side theirs as an alternative formulation of Euclidean geometry {\it per se}.  The intention is quite different.  

In fact, our postulates are the same as those of the Cartesian plane listed above, except that Postulate (2) is replaced by the (ultimately) equivalent,

\begin{enumerate}  \item[(2$'$)]  $\ell$ is continuously differentiable away from the origin, and generates the identity element ``$\mathbf{1}$" of the group action of $\mathbb{R}^\times \hspace{-1mm}\times \mathbb{R}^\times$ on itself, and the trivial homomorphism of $\mathbb{R}^\times \hspace{-1mm}\times \mathbb{R}^\times$ to $\mathbb{R}^\times$, respectively, according to \begin{eqnarray} \label{5/25/25.1}    \mathbf{1} = \ell\cdot\nabla\ell , \\ 1 = \ell\circ\nabla\ell, \label{5/25/25.2}  \end{eqnarray} where the symbol ``$\,\cdot\,$" indicates multiplication of a scalar with a vector, and ``$\,\circ\,$" indicates composition of operations.   The length of a line segment with endpoints $P_1$, $P_2$ is $\ell(P_1-P_2)$.
\end{enumerate}     
 
 Unlike the other formulations of Euclidean geometry, where (at the very least) there is the primitive notion of congruence of line segments, it is not immediately obvious that postulates (1), (2$'$) imply the actual existence of the putative line segment length fulfilling that purpose.   However, with respect to the standard basis, equations (\ref{5/25/25.1}), (\ref{5/25/25.2}) require the non-origin points of the plane to satisfy, \begin{equation}\label{9/6/25.2} (x,y) = \ell(x,y)\nabla\ell(x,y),\end{equation}  \begin{equation} \label{9/6/25.3} \ell(\nabla\ell(x,y)) = 1.\end{equation} On any domain where $\nabla\ell$ {\it might} exist, i.e., where $\nabla\ell = (\frac{\partial\ell}{\partial x},\frac{\partial\ell}{\partial y})$, one can insert the right-hand-side of the latter equation into the above displayed equations, and then easily show that $\ell(x,y)=\sqrt{x^2+y^2}$ (by performing one-dimensional indefinite integration with respect to the  equality of the first components on the left-hand-side and right-hand-side of the resulting (\ref{9/6/25.2}), then independently repeating that procedure with respect to equality of the second components, and finally reconciling the results with each other and (\ref{9/6/25.3})).  One can then immediately verify that $\sqrt{x^2+y^2}$ {\it is} continuously differentiable away from the origin (uniquely so, as a nonnegative function satisfying (\ref{5/25/25.1}), (\ref{5/25/25.2})), and is thus truly the length function.  
 
 So, as opposed to the Cartesian plane formulation that uses Postulate (2) rather than Postulate (2$'$), in this formulation the Pythagorean Theorem actually has to be proved (rather than appropriated from synthetic geometry), and the proof simultaneously assures existence of the length function.    
  Not surprisingly, because the underlying postulates are so different, this proof of the Pythagorean Theorem is quite different from those found in the large compendium \cite{loomis:1940}, the more recent survey \cite{maor:2007}, or the relevant {\it Wikipedia} page \cite{wiki:2025}.  But it has something in common with various of the proofs contained in those sources: it is very simple (assuming one knows Calculus).  
  
   A cynic might (correctly) assert that the form of Postulate (2$'$) is due to prior knowledge of the Pythagorean Theorem.  On the other hand, we would like to think that a superior intelligence would have independently forwarded its equations for the Euclidean norm as an {\it a priori} embodiment of the elegance and simplicity that must accompany the ``most basic" of all norms.
  
 Regardless, the formulation of the Cartesian plane resulting from the four postulates of Birkhoff is equivalent to the alternative formulation resulting from postulates (1) and (2$'$) above.  
 As already indicated, we don't mean to suggest that the latter is anything other than a trivial re-formulation as far as Euclidean geometry itself is concerned.  Instead, the justification for Postulate (2$'$) is that it explicitly forwards a paradigm, i.e., the general format of (\ref{5/25/25.1}), (\ref{5/25/25.2}), to be used in the development of a formalism leading to a new class of algebra isomorphism invariants and accompanying geometries for algebras over the field $\mathbb{R}$.

Some basic results are derived in the next section.  Our program is motivated differently and more fully developed in \cite{fgreensite:2025}.

\bigskip

 {\it Conventions.} Use of the word ``algebra" always refers to a unital associative algebra with vector space of elements $\mathbb{R}^n$, $n\ge 1$.  The standard topology on $\mathbb{R}^n$ is assumed, as well as the standard basis with elements taken to be column vectors.  Superscript $^T$ indicates dual vector or matrix transpose as applicable.  In the sequel, ``$\,\cdot\,$" is the usual dot product, and ``$\mathbf{1}$" is the multiplicative identity of an algebra.  

 \section{The geometries associated with an algebra} 

Equations (\ref{5/25/25.1}), (\ref{5/25/25.2}) are also a prescription for the Euclidean norm in any number of dimensions, \mbox{$\ell:\mathbb{R}^n\rightarrow\mathbb{R}$}, which indeed continues to satisfy,
 \begin{eqnarray}\label{8/5/22.1} s = \ell(s) \nabla \ell(s),\\  \ell(\nabla \ell(s)) = 1, \label{8/5/22.2}\end{eqnarray} for $s\in\mathbb{R}^n$, $s\ne 0$. 
          This general decompositional format, where the non-origin points are expressed by the product of their norm with a point on the unit sphere dependent on the norm's gradient, has relevance well beyond Euclidean space.   To begin with, it applies generally to real quadratic spaces.   Thus, if $L$ is the scalar product matrix associated with a quadratic space where the associated quadratic form is positive-definite or indefinite but nondegenerate, it follows that
  \begin{eqnarray}\label{5/22/25.1} s = \ell_{\rm qs}(s) L^{-1}[\nabla \ell_{\rm qs}(s)],\\ \ell(L^{-1}[\nabla \ell_{\rm qs}(s)]) = 1, \label{5/22/25.11}\end{eqnarray} where the above equations pertain to a domain in $\mathbb{R}^n$ on which the quadratic form is strictly positive and written as $\ell_{\rm qs}^2(s)$, and $\ell_{\rm qs}(s)$ is the positive square-root.  To see this, begin with the quadratic form defined by $ \ell_{\rm qs}^2(s) \equiv s^TLs$.  Take the gradient of both sides of the latter equation and then apply $L^{-1}$ to both sides of the resulting equation to obtain (\ref{5/22/25.1}).  
  Since a quadratic form is a degree-2 positive homogeneous function, it follows that $\ell_{\rm qs}(s)$ is degree-1 positive homogeneous.  Now apply $\ell_{\rm qs}$ to both sides of  (\ref{5/22/25.1}) and invoke the homogeneity of $\ell_{\rm qs}$ to obtain (\ref{5/22/25.11}).  In reference to (\ref{8/5/22.1}), (\ref{8/5/22.2}) and (by analogy) referring to $\ell_{\rm qs}$ as the norm (even though it may now fail the triangle inequality), the points of the relevant domain are again expressed by the product of their norm with a point on the unit sphere dependent on the norm's gradient (where here and in the sequel, ``unit sphere" is always defined with respect to the ``norm" in question).
    
  Although we know that a scalar product matrix is symmetric, the decomposition (\ref{5/22/25.1}), (\ref{5/22/25.11}) actually requires that $L$ be a symmetric matrix.  Indeed, (\ref{5/22/25.1}) implies, \begin{equation}\label{5/26/25.1} Ls = \nabla\left(\frac{1}{2}\ell_{\rm qs}^2(s)\right),\end{equation} so $Ls$ must be a gradient.  Hence, $L$ must satisfy   
  \begin{equation}\label{5/22/25.2} d\big((Ls)\cdot\mathbf{d}s\big) = d\big(s^T L \mathbf{d}s\big) = 0,\end{equation} where $d$ applied to a parenthetical expression indicates exterior derivative, and $\mathbf{d}s\equiv (dx_1,\dots,dx_n)^T$.  We utilize expressions $(Ls)\cdot\mathbf{d}s$ and $s^T L \mathbf{d}s$ to indicate the expected formal symbol manipulations as convenient shorthand for the dual of the vector field $Ls$.  It is easily verified that (\ref{5/22/25.2}) holds if and only if $L$ is symmetric.
  
Now let us suppose the points $s\in\mathbb{R}^n$ are the elements of an algebra and we want a revealing decomposition of the elements along the lines of what is possible for the members of a quadratic space via (\ref{5/22/25.1}), (\ref{5/22/25.11}).  An important distinguishing feature of an algebra, as opposed to only a vector space, is that there are implied nonlinear functions on the elements of the space mapping into the space, e.g., specified by $s\mapsto s^{-1}$, $s\mapsto s^2$, etc.  Let us denote such a function as $f(s)$.  Based on the above, it may be expected that enroute to finding a decomposition similar to (\ref{5/22/25.1}), (\ref{5/22/25.11}) that is influenced by $f$ we will have to find a symmetric matrix $L$ such that \mbox{$d((L[f(s)])\cdot\mathbf{d}s)=0$}.  For definiteness, we can set $f(s)=s^{-1}$, where the domain of $f$ is some neighborhood of the algebra's multiplicative identity $\mathbf{1}$ that is composed only of units (in fact, it might be sufficient to consider {\it only} this function, as indicated in Section 8 of \cite{fgreensite:2025}).  Then, instead of (\ref{5/22/25.2}), we will be dealing with,
  \begin{equation}\label{5/22/25.3} d\left(\left(Ls^{-1}\right)\cdot\mathbf{d}s\right) = d\left((s^{-1})^T L \mathbf{d}s\right) = 0.\end{equation} A solution $L$ to the above equation will be called an uncurling metric, because it annihilates the curl of the vector field $s^{-1}$.   
  
  Since $Ls^{-1}$ is a gradient involving $s^{-1}$ according to (\ref{5/22/25.3}), we can think of it as deriving from a logarithm.  Thus, it makes sense to introduce a function $\ell(s)$ satisfying $\ell(\mathbf{1})=1$ and,
  \begin{equation}\label{11/7/20.1} Ls^{-1}  = \nabla\left[\|{\bf 1}\|^2\log  \ell(s)\right] =   \left(\frac{\|{\bf 1}\|^2}{ \ell(s)}\right)\nabla \ell(s),\end{equation} where $\|{\bf 1}\|^2 \equiv {\bf 1}\cdot {\bf 1}$.  This already seems to have a kinship with (\ref{5/22/25.1}) when the latter is written as $Ls = \ell(s)\nabla\ell(s)$.
  But things become a lot nicer if we consider those uncurling metrics having the additional property whereby elements that are inverses with respect to the algebra's product also behave inversely with respect to the scalar product represented by the uncurling metric.  We refer to the latter as ``normalized" uncurling metrics, denoted by $L_{\shortparallel}$.  So, by definition, a normalized uncurling metric satisfies (\ref{5/22/25.3}) as well as the additional constraint, 
      \begin{equation} \label{7/7/21.1} s^TL_{\shortparallel} s^{-1} =\|{\bf 1}\|^2,\end{equation} in an approriate neighborhood of $\mathbf{1}$. 
 We then specify that $\ell_{\rm sp}(s)$ is a ``special" version of $\ell(s)$ that results when $L$ on the left-hand-side of (\ref{11/7/20.1}) is a normalized uncurling metric $L_{\shortparallel}$.
  
  This results in a very satisfying situation since according to Proposition \ref{9/24/23.1} below, $\ell_{\rm sp}(s^{-1}) = \left(\ell_{\rm sp}(s)\right)^{-1} = \frac{1}{\ell_{\rm sp}(s)}$, and $\ell_{\rm sp}(s)$ is a degree-1 positive homogeneous function.  Hence, using both those properties, and assuming  $L_{\shortparallel}$ is nonsingular, (\ref{11/7/20.1}) then implies
   \begin{eqnarray}\label{5/22/24.5}  s = \ell_{\rm sp}(s)\Big(\|\mathbf{1}\|^2 L_{\shortparallel}^{-1}[\nabla\ell_{\rm sp}(s^{-1})]\Big),\\
\label{5/22/24.6}\ell_{\rm sp}\Big(\|\mathbf{1}\|^2 L_{\shortparallel}^{-1}[\nabla\ell_{\rm sp}(s^{-1})]\Big) = 1,\,\end{eqnarray}
where $\nabla\ell_{\rm sp}(s^{-1})$ means that $\nabla\ell_{\rm sp}$ is evaluated at $s^{-1}$.  In other words, with respect to a neighborhood of $\mathbf{1}$ composed only of units, \begin{quote} There is a potentially large and interesting collection of ways to express an \mbox{algebra} element in a format generally consistent with the paradigm forwarded by the equations of Postulate (2$'$). \end{quote}

The similarity of (\ref{5/22/24.5}), (\ref{5/22/24.6}), with  (\ref{5/22/25.1}), (\ref{5/22/25.11}), is one reason why we have proposed (\ref{11/7/20.1}).  Comparing (\ref{5/26/25.1}) (which pertains to quadratic spaces) with the first equality in (\ref{11/7/20.1}), and further comparing (\ref{5/22/25.1}), (\ref{5/22/25.11}) (which again pertain to quadratic spaces) with  (\ref{5/22/24.5}), (\ref{5/22/24.6}), we are motivated to refer to the algebra-generated spaces described by the latter two equations as ``logarithmic spaces".

Specifically, from (\ref{5/22/25.1}), {\it in a quadratic space the square of the norm} $\ell_{\rm qs}(s)$ is obtained from path-independent integration of $t$ with respect to the metric (symmetric tensor) $L_{\rm sym}$ according to, \begin{equation}\label{8/19/25.1} \big(\ell_{\rm qs}(s)\big)^2 = 2\int_0^s (L_{\rm sym}t)\cdot \mathbf{d}t =  2\int_0^s t^T L_{\rm sym} \mathbf{d}t.\end{equation}
From (\ref{5/22/25.1}), (\ref{5/22/25.11}), we are interested in the unit sphere as the locus of points satisfying $\ell_{\rm qs}(s)=1$ since (as typified with Eucldean space or Minkowski space) the angle associated with two vectors $s_1$, $s_2$ is defined as the arclength of the geodesic path on the unit sphere connecting $\frac{s_1}{\ell_{\rm qs}(s_1)}$, $\frac{s_2}{\ell_{\rm qs}(s_2)}$ where, as with the integration on the right-hand-side of (\ref{8/19/25.1}) that is taken with respect to the metric $L_{\rm sym}$\,, arclength and geodesic are also defined with respect to the metric $L_{\rm sym}$\,. 

Analogously, from (\ref{11/7/20.1}), {\it in a logarithmic space the logarithm of the norm} $\ell_{\rm sp}(s)$ is obtained from path-independent integration of $t^{-1}$ with respect to the uncurling metric $L_{\shortparallel}$ \mbox{according to},
\begin{equation}\label{8/19/25.2} \log\ell_{\rm sp}(s) =  \frac{1}{\|{\bf 1}\|^2} \int_{\bf 1}^s \left(L_{\shortparallel}t^{-1}\right) \cdot \mathbf{d}t =
\frac{1}{\|{\bf 1}\|^2} \int_{\bf 1}^s \left( t^{-1}\right)^T L_{\shortparallel} \mathbf{d}t .\end{equation} 
 From
(\ref{5/22/24.5}), (\ref{5/22/24.6}), the unit sphere is now the locus of points satisfying $\ell_{\rm sp}(s) = 1$.    The angle between two vectors $s_1$, $s_2$ is defined in the same manner as accomplished in a quadratic space, i.e., as the arclength of the geodesic path on the unit sphere connecting $\frac{s_1}{\ell_{\rm sp}(s_1)}$, $\frac{s_2}{\ell_{\rm sp}(s_2)}$ where, as with the integration on the right-hand-side of (\ref{8/19/25.2}) that is taken with respect to the uncurling metric $L_\shortparallel$\,, arclength and geodesic are also defined with respect to the uncurling \mbox{metric $L_{\shortparallel}$}\,.

Thus, as with the distinct geometries arising from different signatures of scalar products defining quadratic spaces, distinct geometries associated with an algebra (in the sense of assigned ``norms" and ``angles") will result from different classes of solutions $L_{\shortparallel}$ satisfying (\ref{5/22/25.3}) and (\ref{7/7/21.1}) - defining the logarithmic spaces.

\begin{definition}\label{9/6/25.1} With respect to algebra $A$, the real vector space of uncurling metrics is the algebra's {\it anti-rotor}, $\mathfrak{u}_A$.  It may also be denoted by $\mathfrak{u}$ if no confusion arises.    \end{definition}

  \begin{theorem} \label{6/6/23.1} If $K:A_1\rightarrow A_2$ is an algebra isomorphism, then \begin{equation}\label{5/23/23.1} \mathfrak{u}_{A_1} = K^T \mathfrak{u}_{A_2} K, \end{equation} 
   where $K^T \mathfrak{u}_{A_2} K \equiv \{K^T L K:L\in \mathfrak{u}_{A_2}\}$.
   \end{theorem}

\begin{proof} By hypothesis, a nonsingular linear change of variable expression for the elements of $A_2$ in terms of the elements of $A_1$ can be written as $s_2= Ks_1$. \mbox{For $L_2\in \mathfrak{u}_{A_2}$} we obtain, \begin{equation} 0 = d\left((s_2^{-1})^T L_2 \mathbf{d}s_2\right) =  d\left(((Ks_1)^{-1})^T L_2 \mathbf{d}(Ks_1)\right) = d\left((s_1^{-1})^T \big[K^T L_2K\big]\mathbf{d}s_1\right), \label{9/21/23.1}  \end{equation} where we have used (\ref{5/22/25.3}), the commutativity of a differential with a linear transformation, and $(Ks_1)^{-1} = K(s_1^{-1})$ (since $K$ is an isomorphism).   It follows from (\ref{9/21/23.1}) that every uncurling metric  of algebra $A_2$ is associated with an uncurling metric of $A_1$, i.e., from the last equality in (\ref{9/21/23.1}) we have  $K^TL_2K\in\mathfrak{u}_{A_1}$.  The isomorphism hypothesis allows us to apply an analogous argument employing $K^{-1}$, indicating that for $L_1\in\mathfrak{u}_{A_1}$ we have the associated $(K^{-1})^T L_1 K^{-1} \in \mathfrak{u}_{A_2}$.    Furthermore, $K^T L_2K\in\mathfrak{u}_{A_1}$ implies $(K^{-1})^T [K^T L_2 K]K^{-1} = (K^T)^{-1} [K^T L_2 K]K^{-1} = L_2$.  Similarly, $(K^{-1})^T L_1K^{-1}\in\mathfrak{u}_{A_2}$ is such that $K^T[(K^{-1})^T L_1K^{-1}] K = L_1$.  The theorem follows from the implied bijection. 
\end{proof}

\begin{corollary}\label{7/31/25.1} If $A_1$ and $A_2$ are isomorphic algebras, then $\mathfrak{u}_{A_1}$ and $\mathfrak{u}_{A_2}$ are the same up to a change \mbox{of basis}. \end{corollary} 

It is seen that (\ref{5/23/23.1}) allows us to immediately identify various algebra isomorphism invariants resulting from the anti-rotor. For example,
            
  \begin{corollary}\label{9/3/25.1}  The following are algebra isomorphism invariants. \begin{itemize} \item The dimension of $\mathfrak{u}$.   \item The largest rank  of a member of $\mathfrak{u}$.  \item The smallest rank of a nonzero member of $\mathfrak{u}$. \end{itemize} \end{corollary}
  
  It is shown in \cite{fgreensite:2025} that computation of $\mathfrak{u}$ is straightforward, and its dimension is always at least 1 but ``usually" greater than 1.

Now let us investigate the question of algebra-related geometries a bit further.

   \begin{definition}\label{9/8/25.1} Let $L$ be an uncurling metric for some algebra $A$.  For $s$ in an open ball centered at the multiplicative identity $\mathbf{1}$ and devoid of non-units, with $\|{\bf 1}\|^2 \equiv \mathbf{1}\cdot\mathbf{1}$, the associated {\it unital norm} $ \ell(s)$ is the function 
  \begin{equation}\label{11/7/20.2} \ell(s) \equiv \exp\left(\frac{1}{\|{\bf 1}\|^2}\int_{\bf 1}^s \left(Lt^{-1}\right)\cdot \mathbf{d}t\right).\end{equation} 
   The set of unital norms associated with an algebra, which is an abelian group $\mathfrak{u}^d_A$ with respect to the usual point-wise operation of multiplication of functions, is called the {\it dual of the anti-rotor}. $\mathfrak{u}^d_A$ may be denoted as $\mathfrak{u}^d$ if no confusion arises.  
   \end{definition} 
   
    Based on the mapping from $\mathfrak{u}$ to $\mathfrak{u}^d$ defined by (\ref{11/7/20.2}),  we can write \begin{equation}  \mathfrak{u}^d = \exp\left(\int_{\bf 1}^s [\mathfrak{u}t^{-1}]\cdot \mathbf{d}t\right). \label{6/7/23.1}\end{equation}  In the above, the notation $\mathfrak{u}t^{-1}$ means $\{Lt^{-1}: L\in \mathfrak{u}\}$.  Accordingly, the right-hand-side of (\ref{6/7/23.1}) is understood to be a collection of exponentiated integrals.   Note that in (\ref{6/7/23.1}) we do not need the factor $\frac{1}{\|{\bf 1}\|^2}$, since $\frac{1}{\|{\bf 1}\|^2}\mathfrak{u}=\mathfrak{u}$.  The entities $\mathfrak{u}$ and $\mathfrak{u}^d$ are considered ``dual" based on the invertible transform relationship (\ref{11/7/20.2}).

Because of the evident analogy with quadratic space geometries, we earlier highlighted the logarithmic space geometries resulting from normalized uncurling metrics $L_\shortparallel$.  It is shown in \cite{fgreensite:2025} that a normalized uncurling metric always exists, and the space of normalized uncurling metrics ``usually" has positive dimension.  Furthermore, it is also shown in \cite{fgreensite:2025} that there is always a normalized uncurling metric associated with a special unital norm closely related to the ``usual" algebra norm (the usual norm being the determinant of an element's image under the left regular representation mapping).  However, any uncurling metric $L$ results in ``lengths" and ``angles" in the same way as the latter are obtained with $L_\shortparallel$.
The following indicates that the resulting collection of geometries persists under an algebra isomorphism.

  \begin{corollary}\label{6/6/23.2} 
If $A_1$ and $A_2$ are isomorphic algebras, then $\mathfrak{u}_{A_1}^d$ and  $\mathfrak{u}_{A_2}^d$ are the same up to a change of basis.
 \end{corollary}

\begin{proof} Let $K:A_1\rightarrow A_2$ be the algebra isomorphism, where for $s_1\in A_1$ we define $s_2\equiv Ks_1$. 
Equation (\ref{5/23/23.1}) implies,
\begin{eqnarray}\nonumber \int_{{\bf 1}_{A_1}}^{s_1} \left(\mathfrak{u}_{A_1}t_1^{-1}\right)\cdot \mathbf{d}t_1 &=& \int_{{\bf 1}_{A_1}}^{s_1} \left( K^T\mathfrak{u}_{A_2}Kt_1^{-1}\right)\cdot \mathbf{d}t_1 \\ \nonumber &=& 
\int_{{\bf 1}_{A_2}}^{s_2} \left( K^T\mathfrak{u}_{A_2}K (K^{-1}t_2)^{-1}\right)\cdot \mathbf{d}\left(K^{-1}t_2\right) \\ &=& \int_{{\bf 1}_{A_2}}^{s_2} \left(K^T\mathfrak{u}_{A_2}t_2^{-1}\right)\cdot \left(K^{-1}\mathbf{d}t_2\right) = \int_{{\bf 1}_{A_2}}^{s_2} \left(\mathfrak{u}_{A_2}t_2^{-1}\right)\cdot \mathbf{d}t_2, \label{9/6/24.2}
\end{eqnarray} 
where we have used the isomorphism assumption (which implies $(K^{-1}t_2)^{-1} = K^{-1}t_2^{-1}$), the commutativity of a linear transformation and a differential, and exploitation of an adjoint in the context of an inner product.  

Recalling (\ref{6/7/23.1}), the corollary follows.
\end{proof}

 Finally, we provide a proof of the feature used to derive (\ref{5/22/24.5}), (\ref{5/22/24.6}).

  \begin{proposition}\label{9/24/23.1} $ \ell_{\rm sp}(s)$ is a degree-1 positive homogeneous function, and
\begin{equation}\label{5/23/23.3}  \ell_{\rm sp}(s^{-1}) = \big(\ell_{\rm sp}(s)\big)^{-1} = \frac{1}{\ell_{\rm sp}(s)}, \end{equation} when $s,s^{-1}$ are in an open ball centered at $\mathbf{1}$ devoid of non-units. \end{proposition}

\begin{proof}  The requisite open ball in the theorem statement always exists, as can be easily demonstrated by considering the left regular representation of the algebra.  Let $\eta$ be this neighborhood of $\mathbf{1}$.

Given $\alpha>0$ such that $\frac{1}{\alpha}s$ and $\alpha s$ are both in $\eta$, the given specifications $\ell_{\rm sp}(\mathbf{1})=1$ and (\ref{11/7/20.1}) imply,
\begin{eqnarray}\nonumber  \log\ell_{\rm sp}(\alpha s) &=&  \frac{1}{\|{\bf 1}\|^2}\int_{\bf 1}^{\alpha s} \left(L_{\shortparallel}t^{-1}\right)\cdot \mathbf{d}t = \frac{1}{\|{\bf 1}\|^2}\int_{\frac{\mathbf{1}}{\alpha}\mathbf{1}}^{s} \left(L_{\shortparallel}(\alpha w)^{-1}\right)\cdot \mathbf{d}(\alpha w)  \\ \nonumber &=&  \frac{1}{\|{\bf 1}\|^2}\int_{\mathbf{1}}^{s} \left(L_{\shortparallel}w^{-1}\right)\cdot \mathbf{d}w + \frac{1}{\|{\bf 1}\|^2}\int_{\frac{\mathbf{1}}{\alpha} \mathbf{1}}^{\mathbf{1}} \left(L_{\shortparallel}w^{-1}\right)\cdot \mathbf{d}w  \\ &=& \log \ell_{\rm sp}(s) + \log\alpha, \label{4/2/23.1}
\end{eqnarray} where the integral from $\frac{\mathbf{1}}{\alpha} \mathbf{1}$ to $\mathbf{1}$ can be easily evaluated along the line segment with those endpoints, and we have used (\ref{7/7/21.1}) in evaluation of that integral to  obtain the term $\log \alpha$ on the right-hand-side of the final equality.  Thus, we have $\ell_{\rm sp}(\alpha s) = \alpha\ell_{\rm sp}(s)$, degree-1 positive homogeneity.

Let $\eta^{-1}$ be the set consisting of the multiplicative inverses of all members of $\eta$.  It is easily shown that $\eta^{-1}$ is an open neighborhood of $\mathbf{1}$, and $\eta$, $\eta^{-1}$ are diffeomorphic via the multiplicative inversion operation.  
We now consider a smooth path $\mathcal{P}$ from $\mathbf{1}$ to $s$ which is inside $\eta$.  Let $\mathcal{P}^{-1}$ be the point set consisting of the multiplicative inverses of the members of $\mathcal{P}$.  Clearly, $\mathcal{P}^{-1}$ is a continuous path contained in $\eta^{-1}$. 
 We then have, 
\begin{eqnarray}  \log \ell_{\rm sp}(s^{-1})   &=&   \frac{1}{\|{\bf 1}\|^2} \int_{\bf 1}^{s^{-1}} \left(L_{\shortparallel}t^{-1}\right) \cdot \mathbf{d}t =  \frac{1}{\|{\bf 1}\|^2} \int_{\bf 1}^s \left(L_{\shortparallel}y\right)\cdot \mathbf{d}\left(y^{-1}\right) \nonumber \\ &=&  \frac{1}{\|{\bf 1}\|^2}\left(\left(y^{-1}\cdot \left(L_{\shortparallel}y\right)\right)\bigg|_{\bf 1}^s - \int_{\bf 1}^s y^{-1}\cdot \mathbf{d}\left(L_{\shortparallel}y\right)\right)   \nonumber \\ &=& -\frac{1}{\|{\bf 1}\|^2} \int_{\bf 1}^s y^{-1}\cdot (L_{\shortparallel}\mathbf{d}y) = -\frac{1}{\|{\bf 1}\|^2} \int_{\bf 1}^s \left(L_{\shortparallel}y^{-1}\right) \cdot \mathbf{d}y  = -\log \ell_{\rm sp}(s),\label{4/2/23.5} \end{eqnarray} where we have used the change of variable $y=t^{-1}$, commutativity of a linear transformation and a differential, the property that $L_{\shortparallel}$ is a real symmetric matrix (i.e., self-adjoint), and  (\ref{7/7/21.1}).  Equation (\ref{5/23/23.3}) then follows.
\end{proof}

   \end{document}